\documentclass{article}

\newcommand{\R}{I\!\! R}
\newcommand{\C}{I\!\! C}
\newcommand{\Z}{I\!\! Z}
\newcommand{\N}{I\!\! N}

\begin{document}
\vspace*{-40pt}
%\begin{center}
\rm

\noindent {\sl \small Nova Science Publishers, Inc.}

\noindent {\sl \small African Diaspora Journal of Mathematics}

\noindent {\sl \small Volume 00, No. 00, p. 000--000 (2004)}

\noindent {\sl \small ISSN 1539-854X \ Copyright $\copyright$2004}

\noindent {\sl \small www.african-j-math.org}

%\end{center}
\vspace*{16pt}

\vspace{1.2cm}

\def\bfE{\mbox{\boldmath$E$}}
\def\bfG{\mbox{\boldmath$G$}}

\title{Non abelian cohomology: the point of view of gerbed tower  }

\author{Aristide Tsemo\thanks{ Department, of
        Mathematics, St George University of Toronto
        Toronto, ON Canada M5S 3G3.
         ({\tt tsemo58@yahoo.ca}).}}

\maketitle

\begin{abstract}
We define in this paper the notion of gerbed tower. This enables
us to interpret   geometrically  cohomology classes without using
the notion of $n$-category. We use this theory to study sequences
of affine maps between affine manifolds, and the cohomology of
manifolds.
\end{abstract}
{\bf keywords}
%\begin{keywords}
gerbes, non abelian cohomology.
%\end{keywords}

{\bf Classification A.M.S.}
%\begin{AMS}
18D05, 57R20.
%\end{AMS}

\pagestyle{myheadings}
\thispagestyle{plain}
 \markboth{A. Tsemo}{{\small Non abelian cohomology: the point of view of gerbed towers}}

\section{Introduction }

\bigskip
\bigskip
In mathematics, a theory  is defined by  axioms which describe
relations between elements of a set. The purpose of mathematicians
is to classify these elements by defining  structures modelled on
reference objects. In geometry, a structure modelled on  the space
$L$, is defined on a topological set $N$, by a Cech $0$-chain
whose boundary  reflects properties of $L$. For example, a
$n$-differentiable manifold is defined by an atlas $(U_i)_{i\in
I}$, and charts $\phi_i:U_i\rightarrow {\R}^n$, such that
$\phi_j\circ\phi_i^{-1}$ is a differentiable map, here the model
$L$ is ${\R}^n$, and the property reflected is the
differentiability. The manifold is obtained by gluing the sets
$\phi_i(U_i)$ using the $1$-cocycle
$h_{ij}=\phi_j\circ\phi_i^{-1}$.  Alternatively, a structure is
defined by gluing a family of sets $N_i$ using a cocycle $h_{ij}$.
Often the sets $N_i$ are related to the model in the sense that
each of them is endowed with a $L$-structure. The natural problem
to determine if a given topological space $N$ can be endowed with
a structure modelled on $L$, leads to the notion of sheaf of
categories. When the structure exists locally, that is when there
exists  an open cover $(N_i)_{i\in I}$ of $N$, such that each
$N_i$ is endowed with a $L$-structure, the existence of the
$L$-structure on $N$ is equivalent to determine whether  the
cohomology class of a $2$-Cech cocycle is trivial. This has
motivated the definition of a $2$-structure called gerbe, which is
classified in geometry by a $2$-Cech cocycle. The natural problem
which occurs is to provide geometric conditions which insure a
$2$-Cech chain to define a $2$-type structure,..., a $n$-Cech
chain to define a $n$-type structure. On this purpose, one needs
to give a geometric interpretation of Cech cohomology classes.
Unfortunately, the notion of $n$-category needed to define
$n$-structures is not well-understood. The main goal of this paper
is to interpret Cech classes geometrically, by defining the notion
of commutative $n$-gerbed tower. These are sequences of
$2$-categories $F_n\rightarrow F_{n-1}...F_2\rightarrow F_1$,
where $F_1$ is a gerbe defined on a topos $N$. A commutative
$n$-gerbed tower satisfy conditions which allow to attach to it a
family of  $p$-Cech cohomology classes $([f_2],...,[f_{n+1}])$,
where $[f_p]\in H^{p}(N,L_{p-1})$, and $L_p$ are commutative
sheaves defined on $N$. This notion represents geometrically the
connecting morphism in cohomology. More precisely we have:

\medskip

{\bf Theorem 4.2.6.}

{\it Let $u_n=F_n\rightarrow F_{n-1}...F_2\rightarrow F_1$ be a
commutative $n$-gerbed tower defined on a topos to which is
associated the family of cohomology classes
$([f_2],...,[f_{n+1}])$. Suppose that there exists an exact
sequence of sheaves $0\rightarrow L_{n+1}\rightarrow
L'_{n+1}\rightarrow L_n\rightarrow 0$, then the family of
cohomology classes $([f_2],...,[f_{n+2}])$ where $[f_{n+2}]$ is
the image of  $[f_{n+1}]$ by the connecting morphism
$H^{n+1}(N,L_n)\rightarrow H^{n+2}(N,L_{n+1})$, is associated to a
$n+1$-gerbed tower.}

\medskip

 An example of a $n$-gerbed tower appears in the theory of affine
manifolds.
 An affine manifold $(N,\nabla_N)$ is a differentiable manifold
$N$, endowed with a connection $\nabla_N$, whose curvature and
torsion forms vanish identically. We say that the $n$-dimensional
affine manifold $(N,\nabla_N)$, is complete, if and only if it is
the quotient of the affine space ${\R}^n$, by a subgroup
$\Gamma_N$ of $Aff({\R}^n)$ which acts properly and freely on
${\R}^n$. L. Auslander has conjectured that the fundamental group
of a compact and complete affine manifold is polycyclic. In [26]
we have conjectured that a finite Galois cover of a compact and
complete affine manifold is the domain of a non trivial affine
map. This leads to the following problem: classify sequences
$(N_n,\nabla_{N_n})\rightarrow...\rightarrow (N_1,\nabla_{N_1})$
where each map $f_i:(N_{i+1},\nabla_{N_{i+1}})\rightarrow
(N_i,\nabla_{N_i})$ is an affine fibration whose domain is compact
and complete: This means that $f_i$ is a surjective map and each
affine manifold $(N_i,\nabla_{N_i})$ is a compact and complete
affine manifold. The classification of affine fibrations has been
done using gerbe theory (see [28]). It is normal to think that
 composition sequences of affine manifolds are related to
$n$-gerbes.  We define a $n$-gerbed tower which appears naturally
in this context.

 Characteristic classes are used in mathematics to study
many objects, for example, Witten has used  characteristic classes
to study the Jones polynomial. This shows the necessity to give a
geometric interpretation of characteristic classes. On this
purpose, we have to  interpret geometrically  the integral
cohomology of a differentiable manifold $N$. The theory of
Kostant-Weil gives an interpretation of the group $H^2(N,{\Z})$:
It is the set of equivalence classes of complex line bundles over
$N$. In [5], is given an interpretation of $H^3(N,{\Z})$ in terms
of equivalence classes of Dixmier-Douady groupoids. Our theory
enables us to interpret a subgroup of $H^{n+2}(N,{\Z})$ as the set
of equivalence classes of a family of  $n$-gerbed towers.

\medskip

This is the plan of our paper:

1. Introduction.

2. The notion of gerbe.

3.Notations

4. The notion of gerbed towers.

5. Spectral sequences and gerbed towers.

6. Applications of gerbed towers to affine geometry.

7. Interpretation of the integral cohomology of a manifold.

8 A definition of a notion of sheaf of $n$-categories.

\bigskip

\bigskip

\section{ The notion of gerbe.}

\bigskip

In this part we present the notion of gerbe studied by Giraud
[11].

\medskip

{\bf Definitions 2.0.1.}

Let $E$ be a category, a {\bf sieve} is a subclass $R$ of the
class of objects $Ob(E)$ of $E$ such that if $f:X\rightarrow Y$ is
a map of $E$, such that $Y\in R$, then $X\in R$.

  Let $f:E'\rightarrow E$ be a functor, and  $R$ a sieve of $E$,
we denote by $R^f$, the sieve defined by  $R^f=\{X\in Ob(E'):
f(X)\in R\}$.

For each object $T$ of $E$, we denote by $E_T$, the category whose
objects are arrows $u:U\rightarrow T$, a morphism of $E_T$ between
$u_1:U_1\rightarrow T$, and $u_2:U_2\rightarrow T$, is a map
$h:U_1\rightarrow U_2$ such that $u_2\circ h=u_1$ $\bullet$

\medskip

{\bf Definition 2.0.2.}

A {\bf topology} on $E$ is defined as follows:  to each object $T$
of $E$, we associate a non empty set  $J(T)$ of  sieves of the
category $E_T$ of $E$, above $T$ such that:

(i) For each map $f:T_1\rightarrow T_2$, and for each element $R$
of $J(T_2)$, $R^f\in J(T_1)$. (The morphism $f$ induces a functor
between $E_{T_1}$ and $E_{T_2}$ abusively denoted $f$).

(ii) The sieve $R$ of $E_T$  is an element of $J(T)$, if for every
map $f:T'\rightarrow T$ of $E$, ${R}^f\in J(T')$ $\bullet$

\medskip

 A category endowed with a topology is called a site.

\medskip

{\bf Definitions 2.0.3.}

A {\bf sheaf  of sets} $L$ defined on the category $E$ endowed
with the topology $J$, is a contravariant functor $L:E\rightarrow
Set$, where $Set$ is the category of sets, such that for each
object $U$ of $E$, and each element $R$ of $J(U)$, the natural
map:

$$
L(U)\longrightarrow lim(L\mid R)
$$

is bijective, where $(L\mid R)$ is the correspondence defined on
$R$ by $(L\mid R)(f)=L(T)$ for each map $f:T\rightarrow U$ in $R$.

 Let $h:F\rightarrow E$ be a functor, for each object $U$ of $E$,
we  denote by $F_U$ the subcategory    of $F$ defined as follows:
an object $T$ of $F_U$ is an object $T$ of $F$  such that
$h(T)=U$. A map $f:T\rightarrow T'$   between a pair of objects
$T$ and $T'$ of $F_U$, is a map of $F$ such that $h(f)$ is the
identity of $U$. The category $F_U$ is called {\bf the fiber}  of
$U$. For each objects $X$, and $Y$ of $F_U$, we will denote by
$Hom_U(X,Y)$ the set of morphisms of $F_U$ between $X$ and $Y$
$\bullet$

\medskip

{\bf Definitions 2.0.4.}

 Let $h:F\rightarrow E$ be a functor, $m:x\rightarrow y$ a map of
$F$, and $f=h(m):T\rightarrow U$ its projection by $h$. We will
say that $m$ is {\bf cartesian}, or that $m$ is the {\bf inverse
image} of $f$ by $h$, or $x$ is an inverse image of $y$ by $h$, if
for each element $z$ of $F_T$, the map
$$
Hom_T(z,x)\rightarrow Hom_f(z,y)
$$
$$
n\rightarrow mn
$$
is bijective, where $Hom_f(z,y)$ is the set of maps
$g:z\rightarrow y$ such that $h(g)=f$.

 A functor $h:F\rightarrow E$ is a {\bf fibered category} if and
only if each map $f:T\rightarrow U$, has an inverse image, and the
composition of two cartesian maps is a cartesian map.

 We will say that the category is fibered in groupoids, if for
each diagram

$$
x{\buildrel{f}\over{\longrightarrow}}
z{\buildrel{g}\over{\longleftarrow}} y
$$

of $F$ above the diagram of $E$,
$$
U{\buildrel{\phi}\over{\longrightarrow}}
W{\buildrel{\psi}\over{\longleftarrow}} V
$$
 and for each map $m:U\rightarrow V$ such that $\psi m=\phi$, there
 exists a unique map $p:x\rightarrow y$, such that $gp=f$, and
 $h(p)=m$.

 This implies that the inverse image is unique up to
 isomorphism.

 Consider a map $\phi:U\rightarrow V$  of $E$, we can define a functor
 $\phi^*:F_V\rightarrow F_U$, such that for each object $y$ of
 $F_V$, $\phi^*(y)$ is defined as follows: we consider a cartesian map
 $f:x\rightarrow y$ above $\phi$ and set  $\phi^*(y)=x$. Remark
 that although the definition of $\phi^*(y)$ depends of the chosen
 inverse image $f$, the functors $(\phi\psi)^*$ and $\psi^*\phi^*$
 are isomorphic $\bullet$

 \medskip

 {\bf Definitions 2.0.5.}

 A {\bf section} of a fibered category $h:F\rightarrow E$, is a
 correspondence defined on the class of arrows of $E$ as follows:
 to each map $f:U\rightarrow T$, we define a cartesian map:
 $l^f:x_U\rightarrow y_T$ of $F$, whose image by $h$ is $f$ such that:
 $l^{f'f}=l^{f'}\circ l^f$.

\medskip

Consider the diagram
$$
\matrix{F && G\cr \downarrow f && \downarrow g\cr E_1 &
{\buildrel{u}\over{\longrightarrow}} & E_2}
$$
where the functors $f$ and $g$ are fibered categories. We denote
by $Hom_u(F,G)$ the subcategory of $Hom(F,G)$ whose objects are
functors $v:F\rightarrow G$ which verify $gv=uf$. The maps of this
category are morphisms $m:v\rightarrow v'$ such that $gm$ is the
identity morphism of the functor $uf$.

We denote by $Cart_u(F,G)$, the subcategory of $Hom_u(F,G)$ whose
objects are cartesian functors: These are functors which transform
cartesian maps to cartesian maps $\bullet$

\medskip

Let $E$ be a category endowed with a topology $J$, and
$F\rightarrow E$ a fibered category, for each object $U$ of $E$,
and each element $R$ of $J(U)$ ,we consider the canonical functors
$E_U\rightarrow E$, and $R\rightarrow E$. We can define the set of
cartesian functors $Cart_{Id_E}(E_U,F)$ and $Cart_{Id_E}(R,F)$.
There exists a canonical restriction functor
$Cart_{Id_E}(E_U,F)\rightarrow Cart_{Id_E}(R,F)$.

\medskip

{\bf Definition 2.0.6.}

Let $E$ be a category endowed with a topology, a {\bf sheaf of
categories} on $E$, is a fibered category $F\rightarrow E$, such
that for each sieve $R$,
 the cartesian functor $Cart_{Id_E}(E_U,F)\rightarrow Cart_{Id_E}(R,F)$
 defined at the paragraph above is an equivalence of categories
 $\bullet$

\bigskip

{\bf Proposition-Definition 2.0.7.}

{\it Suppose that  $E$ is a topos whose topology is generated by a
contractible covering family $(U_i\rightarrow U)_{i\in I}$, and
$h:F\rightarrow E$ a fibered category in groupoids. For each map
$f:U\rightarrow V$ of $E$, we consider the functor
$r_{U,V}(f):F_V\rightarrow F_U$ defined as follows: For each
object $y$ of $F_V$, $r_{U,V}(f)(y)$ is an object $x$ of $F_U$
such that there exists a cartesian map $n:x\rightarrow y$ such
that $h(n)=f$. Consider the maps $v_1:U_1\rightarrow U_2$, and
$v_2:U_2\rightarrow U_3$ of $E$, the functors
$r_{U_1,U_2}(v_1)\circ r_{U_2,U_3}(v_2)$ and $r_{U_1,U_3}(v_2v_1)$
are isomorphic (see [11]). The functor $h:F\rightarrow E$ is a
sheaf of categories if and only if the correspondence
$U\rightarrow F_U=F(U)$ satisfies the following properties:

(i) Gluing condition for arrows.

Let $U$ be an object of $E$, and $x$, $y$ objects of $F(U)$. The
functor  from $E_U$, endowed with the restriction of the topology
$J$, to the category of sets which associates to an object
$f:V\rightarrow U$ the set $Hom_V(r_{V,U}(f)(x),r_{V,U}(f)(y))$ is
a sheaf of sets.

(ii) Gluing condition for objects.

Consider a covering family $(U_i\rightarrow U)_{i\in I}$ of an
object $U$ of $E$, and for each $U_i$, an object $x_i$ of
$F({U_i})$. Let $t_{ij}:x_j^i\rightarrow x_i^j$, a map between the
respective restrictions of $x_j$ and $x_i$ to $U_i\times_UU_j$
such that on $U_{i_1}\times_UU_{i_2}\times_U U_{i_3}$, the
restrictions of the arrows $t_{i_1i_3}$ and $t_{i_1i_2}t_{i_2i_3}$
are equal.  There exists an object $x$ of $F(U)$ whose restriction
to $F({U_i})$ is $x_i$.

\medskip

If moreover the following properties are verified:

(iii) There exists a covering family $(U_i\rightarrow U)_{i\in I}$
of $E$ such that $F({U_i})$ is not empty,

(iv) For each pair of objects $x$, and $y$ of $F({U_i})$,
$Hom_{U_i}(x,y)$ is not empty (local connectivity),

(v) The elements of $Hom_{U_i}(x,y)$ are invertible. The fibered
category is called a {\bf gerbe}.

(vi) We say that the gerbe is {\bf bounded} by the sheaf $L_F$
defined on $E$, or that $L_F$ is {\bf the band} of the gerbe, if
and only if there exists a sheaf of groups $L_F$ defined on $E$
such that for each object $x$ of $F(U)$ we have an isomorphism:
$$
L_F(U)\rightarrow Hom_U(x,x)
$$
which commutes with restrictions, and with morphisms between
objects$\bullet$}

\bigskip

\subsection{ Classifying cocycle and classification of gerbes.}

\medskip

In this paragraph, we recall  the definition of the classifying
cocycle  of a gerbe defined on the topos $E$ and bounded by the
sheaf  $L$.

\medskip

{\bf Definitions 2.1.1.}

- A gerbe $F\rightarrow E$ is trivial if it has a section.  This
means that $F_E$ is not empty.

- Two gerbes $F\rightarrow E$, and $F'\rightarrow E$ whose band is
$L$, are equivalent if and only if there exists an isomorphism
between the underlying fibered categories which commutes with the
action of $L$. We denote by $H^2(E,L)$ the set of equivalence
classes of gerbes defined on $E$  bounded by $L$ $\bullet$

\medskip

Suppose that the topology of $E$ is defined by the covering family
$(U_i\rightarrow U)_{i\in I}$,  the class of objects of $F_{U_i}$
is not empty, and each objects $x$ and $y$ of $F_{U_i}$ are
isomorphic. Let $(x_i)_{i\in I}$ be a family of objects of $F$,
such that $x_i$ is an object of $F_{U_i}$. There exists a map
$u_{ij}:x_j^i\rightarrow x_i^j$ between the respective
restrictions of $x_j$ and $x_i$ to $F_{U_i\times_UU_j}$. We denote
by $u_{i_1i_2}^{i_3}$ the map between the respective restrictions
of $x_{i_2}$ and $x_{i_1}$ to
$F_{U_{i_1}\times_UU_{i_2}\times_UU_{i_3}}$.

\medskip

{\bf Theorem [11] 2.1.2.}

{\it The family of maps
$c_{i_1i_2i_3}=u^{i_2}_{i_3i_1}u^{i_3}_{i_1i_2}u^{i_1}_{i_2i_3}$
is the classifying $2$-Cech cocycle of the gerbe. If the band $L$
is commutative, then the set of equivalence classes of gerbes over
$E$ whose band is $L$, is one to one with the Cech cohomology
group $H^2(E,L)$.}

\bigskip

\bigskip

\section{ Notations.}

\medskip

Let $U_{i_1},...,U_{i_p}$ be objects  of a topos  $E$, and $C$ a
presheaf of categories defined on $E$. We will denote by
$U_{i_1..i_p}$ the fiber product of $U_{i_1}$,...,$U_{i_p}$ over
the final object. If $e_{i_1}$ is an object of $C(U_{i_1})$,
${e_{i_1}}^{i_2...i_p}$ will be the restriction of $e_{i_1}$ to
$U_{i_1...i_p}$. For a map $h:e\rightarrow e'$ between two objects
of $C(U_{i_1..i_p})$, we denote by $h^{i_{p+1..i_n}}$ the
restriction of $h$ to a morphism between
$e^{i_{p+1}...i_n}\rightarrow {e'}^{i_{p+1}...i_n}$.

\bigskip

\section{ Gerbed tower.}

\bigskip

The purpose of this part is to generalize the notion of gerbe to
the notion of gerbed tower. This notion will allow us to define,
and to represent geometrically higher non abelian  cohomological
classes. In the sequel we assume known the notion of $2$-category
or bicategory defined by Benabou [1]. Recall that a $2$-category
$C$ is defined by a class of objects $Obj(C)$, and for each
objects $x$ and $y$ of $C$, a category $Hom_C(x,y)$ called the
category of morphisms. The objects of $Hom_C(x,y)$ are called
$1$-arrows, and the arrows are called $2$-arrows. There exists a
composition functor:

$$
c(u_1,u_2,u_3):Hom(u_2,u_3)\times Hom(u_1,u_2)\longrightarrow
Hom(u_1,u_3)
$$

For each quadruple $(u_1,u_2,u_3,u_4)$ in $C$, there exists an
isomorphism $c(u_1,u_2,u_3,u_4)$ between the functors

$$
(Hom(u_3,u_4)\times Hom(u_2,u_3))\times
Hom(u_1,u_2)\longrightarrow Hom(u_1,u_4)
$$

and

$$
Hom(u_3,u_4)\times (Hom(u_2,u_3)\times
Hom(u_1,u_2))\longrightarrow Hom(u_1,u_4)
$$

which satisfies more compatibility axioms which can be found in
Benabou. We will suppose that $c(u_1,u_2,u_3,u_4)$ is the identity
on objects. This implies that we can define the category  $C_1$
whose objects are the objects of $C$, and such that for each pair
of objects $x$, $y$ of $C_1$, $Hom_{C_1}(x,y)$ is the set of
objects of $Hom_C(x,y)$.  Let $h:F\rightarrow E$ be a gerbe. We
can define the $2$-category $C(E,F)$ whose objects are objects of
$F$.  Let $x$ and $y$ be a pair of objects of $C(E,F)$, an object
of $Hom_{C(E,F)}(x,y)$ is an arrow between $h(x)$ and $h(y)$. A
$2$-arrow between the objects $x$ and $y$ is a cartesian map
between $x$ and $y$.

\medskip

{\bf Definition 4.0.1.}

  A {\bf bicategory $C$, endowed with a topology} $J$, is a bicategory
  whose objects  are toposes, and for each pair of objects $x$ and $y$ of
  $C$. The set of   $2$-arrows  between $x$ and $y$ is contained in the space
  of continuous maps between $x$ and $y$  $\bullet$

\medskip
\medskip

{\bf Definition 4.0.2.}

A {\bf $n$-gerbed tower} is defined by:

1. A family    $F_n,F_{n-1},... F_2,F_1$  of $2$-categories
respectively endowed with  topologies $J_n$,...,$J_1$, and a
family of $2$-functors $p_l:F_l\rightarrow F_{l-1}$,
$l\in\{2,...,n\}$ which satisfy the following conditions:

\medskip

 2. $F_1$ is a gerbe $p_1:F\rightarrow E$, since we can assume that
 $E$ is a $2$-category such that for each objects $U$ and $V$ of
 $E$, the set of arrows between a pair of elements $f$ and $f'$ of
 $Hom_E(U,V)$ is a singleton, we will often consider the sequence
 of $2$-categories $F_n\rightarrow..F_1 \rightarrow E$. We suppose
 that the $2$-arrows of $F_p$, $p\in \{1,...,n\}$ are invertible.

\medskip

3. Let $U$ be an object of $F_p$, and $l\geq p$ a pair of integers
inferior to $n$. We denote by ${F_{lp}}_U$, the $2$-category whose
class of objects is contained in the class of objects of  $F_l$,
and such that $V$ is an object of ${F_{lp}}_U$ if and only if
$p_{p+1}..p_l(V)=U$. The category of morphisms
$Hom_{{F_{lp}}_U}(X,Y)$ between a pair of objects $X$ and $Y$ is
the subcategory of $Hom_{F_l}(X,Y)$ such that the projections of
$1$-arrows of $Hom_{{F_{lp}}_U}(X,Y)$ by $p_{p+1}..p_l$ is
 $Id_{U}$ of $U$, and the projections of $2$-arrows of
 $Hom_{{F_{lp}}_U}(X,Y)$ by the same functor is the identity of $Id_U$. We
 denote ${F_{l0}}_U$ by ${F_l}_U$. We suppose that for
 each arrow $f:U\rightarrow V$ of $E$, there exists a restriction
 functor $r^l_{U,V}(f):{F_l}_V\rightarrow {F_l}_U$ such that for
 every map $g:V\rightarrow V'$,
 $r^l_{U,V}(f)\circ r^l_{V,V'}(g)=r^l_{U,V'}(gf)$.

\medskip

4.   There exists a family of sheaves $L_1$,...,$L_n$ defined on
$E$. The sheaf $L_{l+1}$ induces a  sheaf ${L_{l+1}}_{U_l} $ on
the object $U_l$ of $F_l$ (recall that $U_l$ is a topos) defined
by its global sections ${L_{l+1}}_{U_l}=L_{l+1}(p_1..p_l(U_l))$.
For each object $U_{l-1}$ of $F_{l-1}$, we suppose that the fiber
of ${F_{ll-1}}_{U_{l-1}}$ is a gerbe defined on the topos
$U_{l-1}$ bounded by  ${L_l}_{U_{l-1}}$.

\medskip

5. Let ${U_1}_l$ and ${U_2}_l$ be a pair of objects of $F_l$, and
$u^2_l:h^1_l\rightarrow {h}^2_l$ a $2$-arrow, between ${U_1}_l$
and ${U_2}_l$, that is an arrow of the category
$Hom_{F_l}({U_1}_l,{U_2}_l)$ between the objects $h^1_l$ and
${h}^2_l$. Recall that $u^2_l$ is a continuous functor between the
topoi ${U_1}_l$ and ${U_2}_l$. We suppose that for every object
${U'_2}_l$ of ${F_{l+1l}}_{{U_2}_l}$ there exists an object
${U'_1}_l$ of ${F_{l+1l}}_{{U_1}_l}$, and a $2$-arrow
${u^2_l}^*:{{U'_1}_l}\rightarrow {{U'_2}_l}$ such that the
following diagram is commutative:

$$
\matrix{{{U'_1}_l}& {\buildrel{{u^2_l}^*}\over{\longrightarrow}}&
{{U'_2}_l}\cr \downarrow p_l && \downarrow p_l\cr {U_1}_l &
{\buildrel{u^2_l}\over{\longrightarrow}} & {U_2}_l}
$$

This implies that for every $2$-arrow $v^2_l:h^2_l\rightarrow
h^3_l$ between the pair of objects ${U_2}_l$ and ${U_3}_l$,  for
every object ${U'_3}_l$ of ${F_{l+1l}}_{{U_3}_l}$, and every
$2$-arrow ${v^2_l}^*$  of $F_{l+1}$ defined by the diagram above,
there exists  an automorphism over the identity $c(u^2_l,v^2_l)$
of ${U'_3}_l$ such that

$$
{v^2_l}^*\circ{u^2_l}^*=c(u^2_l,v^2_l)({v^2_l\circ u^2_l})^*,
$$
since we have supposed that the $2$-arrows are invertible
morphisms of topoi
 $\bullet$

\medskip

{\bf Definition 4.0.3.}

An $\infty$-gerbed tower, is a sequence of functors between
$2$-categories $...F_n\rightarrow F_{n-1}...\rightarrow
F_1\rightarrow E$ such that for each integer $n$, $F_n\rightarrow
F_{n-1}...\rightarrow F_1\rightarrow E$ is a gerbed tower
$\bullet$

\medskip

{\bf Definitions 4.0.4.}

- A morphism $F$ between the gerbed towers $f=F_n\rightarrow..
F_1\rightarrow E$ and $f'=F'_n\rightarrow..\rightarrow
F'_1\rightarrow E$, is defined by a family of $2$-functors
$f_l:F_l\rightarrow F'_l$ such that for each $l$, the following
diagram is commutative:

$$
\matrix{F_l &{\buildrel{f_l}\over{\longrightarrow}} & F'_l\cr
\downarrow p_l && \downarrow p'_l\cr F_{l-1} &
{\buildrel{f_{l-1}}\over{\longrightarrow}} & F'_{l-1}}
$$

and for each object $U_l$ of $F_l$, the induced morphism:
${F_{l+1l}}_{U_l}\rightarrow {F'_{l+1l}}_{f_l(U_l)}$ is a morphism
of gerbes.

- The morphism defined by the family of $2$-functors
$(f_n,..,f_1)$ is an isomorphism, if and only if there exists a
morphism between the gerbed tower $F'_n\rightarrow..\rightarrow
F'_1\rightarrow E$ and $F_n\rightarrow...\rightarrow
F_1\rightarrow E$ defined by the family of $2$-functors
$f'_n,..,f'_1$ such that for each $l$, $f'_l\circ f_l=Id_{F_l}$,
and $f_l\circ f'_l=Id_{F'_l}$.

We say that the gerbed towers $f'$ and $f$ are weakly equivalent
if and only if $f'_l\circ f_l$ is isomorphic to $Id_{F_l}$, and
$f_l\circ f'_l$ is isomorphic to $Id_{F'_l}$. We denote by
$H^{n+1}(E,L_n)$ the set of weakly equivalence classes of gerbed
towers bounded by $(L_1,..,L_n)$. Here $L_n$ is a fixed sheaf
defined on $E$ $\bullet$

\medskip

\subsection{ Non commutative cohomology of groups.}

\medskip

Let $H$ be a group, $V$ a vector space, $Gl(V)$ the group of
linear automorphisms of $V$, and $\rho: H\rightarrow Gl(V)$ a
representation. To define the cohomology groups $H^n(H,V,\rho)$ of
the representation $\rho$, one can consider $EH$ the
$1$-Eilenberg-Maclane space defined by $H$, the representation
$\rho$ defines on $EH$ a flat $V$-bundle whose holonomy is $\rho$.
The cohomology groups $H^n(H,V,\rho)$, are the $n$-cohomology
groups, of the sheaf of locally constant sections of this flat
bundle. This motivates the following definition:

\medskip

{\bf Definition 4.1.1.}

Consider the groups $H$ and $G$, $Aut(G)$ the group of
automorphisms of $G$, and $\rho:H\rightarrow Aut(G)$ a
representation. The representation $\rho$ defines on $EH$ a flat
$G$-bundle $p_{G}$. We denote by $L_{p_{G}}$ the sheaf of locally
constant sections of this bundle. We define $H^{n+1}(H,G,\rho)$ to
be set of weakly equivalence classes of gerbed towers of rank $n$
$F_n\rightarrow..\rightarrow F_1\rightarrow EH$ bounded by a
sequence $(L_1,...,L_{n-1},L_{p_G})$ $\bullet$

\bigskip

\subsection{ The classifying cocycle of a gerbed tower.}

\medskip

Let $f=F_n\rightarrow F_{n-1}\rightarrow... F_1\rightarrow E$ be a
gerbed tower. We will define in this part the classifying cocycle
of $f$. We suppose that the sheaves $L_1,...,L_n$ are commutative,
and there  exists commutative sheaves $L'_1$,...,$L'_n$ defined on
$E$ such that for every objects $U_l$ of $F_l$,  $U_{l+1}$ of
${F_{l+1l}}_{U_l}$, and $h_l$ an object of
$Hom_{F_{l+1}}(U_{l+1},U_{l+1})$,
$Aut({F_{l+1l}}_{h_{l+1}})=L'_{l+1}(p_1..p_l(U_{l}))$ where
$Aut({F_{l+1l}}_{h_{l+1}})$  is the group of automorphisms of the
object ${h_{l+1}}$ whose image by $p_{l+1}$ are elements of
$L_l(p_1...p_l(U_{l}))$.

\medskip

Suppose that the topology of $E$ is defined by the covering family
$(U_i\rightarrow U)_{i\in I}$ such that for each $i$,
${F_1}_{U_i}$ is  not empty, and its objects are isomorphic. Let
$u_i$ be an object of ${F_1}_{U_i}$, and $v_{i_1i_2}$ an arrow
between $u^{i_1}_{i_2}$ and $u^{i_2}_{i_1}$. The family of arrows
$c_{i_1i_2i_3}={v_{i_3i_1}}^{i_2}{v_{i_1i_2}}^{i_3}{v_{i_2i_3}}^{i_1}$
is the classifying cocycle of the gerbe $F_1\rightarrow E$.

If we identify $F_1\rightarrow E$ with a $2$-category, then
$c_{i_1i_2i_3}$ is a $2$-arrow. Let $u_{i_1i_2i_3}$ be an object
of the fiber ${F_{21}}_{u^{i_1i_2}_{i_3}}$. The property 5 of the
definition of gerbed towers implies the existence of a $2$-arrow
${c_{i_1i_2i_3}}^*$ of $Hom_{F_2}(u_{i_1i_2i_3},u_{i_1i_2i_3})$
over $c_{i_1i_2i_3}$.

we can define the automorphism $c_{i_1i_2i_3i_4}$  of the object
${u_{i_1i_2i_3}}^{i_4}$ by:

$$
c_{i_1i_2i_3i_4}={{c_{i_2i_3i_4}}^*}^{i_1}-{{c_{i_1i_3i_4}}^*}^{i_2}+{{c_{i_1i_2i_4}}^*}^{i_3}
-{{c_{i_1i_2i_3}}^*}^{i_4}
$$

each member of the right part of the previous equality can be
supposed to be a morphism of the same object of
$Hom_{F_2}(u_{i_1i_2i_3},u_{i_1i_2i_3})$. The property 5 of the
definition of gerbed towers implies that $c_{i_1i_2i_3i_4}$ is an
element of
$L_2(U_{i_1}\times_EU_{i_2}\times_EU_{i_3}\times_EU_{i_4})$.

\medskip

{\bf Proposition 4.2.1.}

{\it The family $c_{i_1i_2i_3i_4}$ that we have just defined is a
$2$-Cech cocycle.}

\medskip

{\bf Proof.}

The Cech boundary of $c_{i_1i_2i_3i_4}$ is:

$$
\partial(c_{i_1..i_4})=\sum_{p=1}^{p=5}(-1)^{j}c_{i_1..\hat i_p..i_5}
$$

$$
=\sum_{p=1}^{p=5}(-1)^p(\sum_{c=1}^{c=p-1}(-1)^cc_{i_1..\hat
i_c..\hat i_p..}^*+\sum_{c=p+1}^{c=5}(-1)^{c+1}c_{i_1..\hat
i_p..\hat i_c..i_5}^*)=0.
$$

The last sum is zero because the sheaf $L'_1$ is
commutative$\bullet$

\medskip

Suppose defined the classifying cocycles
$c_{i_1i_2i_3}$,...,$c_{i_1..i_{l+2}}$, $l\geq 2$ of the gerbed
tower $F_l\rightarrow...\rightarrow F_1\rightarrow F$. The arrow
$c_{i_1...i_{l+2}}$ is a $2$-arrow of
${u_{i_1..i_{l+1}}}^{i_{l+2}}$. Let $u_{i_1...i_{l+2}}$ be an
object of ${F_{l+1l}}_{{u_{i_1...i_{l+1}}}^{i_{l+2}}}$. The
property 5 implies the existence of an automorphism
${c_{i_1..i_{l+2}}}^*$ of a $1$-arrow of the object
$u_{i_1..i_{l+2}}$ over $c_{i_1..i_{l+2}}$. We can define:

$$
c_{i_1..i_{l+3}}=\sum_{p=1}^{l+3}(-1)^p{c_{i_1..\hat
i_p..i_{l+3}}}^*
$$

We can apply the property 5 to identify $c_{i_1..i_{l+3}}$ to a
$2$-arrow of  ${u_{i_1..i_{l+2}}}^{i_{l+3}}$.

\medskip

{\bf Proposition 4.2.2.}

{\it The family of arrows $c_{i_1..i_{l+3}}$ that we have just
defined is a $l+1$-Cech cocycle.}

\medskip

{\bf Proof.}

The Cech boundary of $c_{i_1..i_{l+3}}$ is:

$$
\partial(c_{i_1..i_{l+3}})=\sum_{p=1}^{p=l+4}(-1)^{p}c_{i_1..\hat i_p..i_{l+4}}
$$

$$
=\sum_{p=1}^{p=l+4}(-1)^p(\sum_{c=1}^{c=p-1}(-1)^cc_{i_1..\hat
i_c..\hat i_p..}^*+\sum_{c=p+1}^{c=l+4}(-1)^{c+1}c_{i_1..\hat
i_p..\hat i_c..i_{l+4}}^*)=0.
$$

The last sum is zero because the sheaf $L'_{l+1}$ is commutative
$\bullet$

\medskip

{\bf Proposition 4.2.3.}

{\it  The cohomology class $[c_{l+2}]$, is the image of
$[c_{l+1}]$ by the connecting morphism $H^{l+1}(E,L_l)\rightarrow
H^{l+2}(E,L_{l+1})$ of the exact sequence $0\rightarrow
L_{l+1}\rightarrow L'_{l+1}\rightarrow L_l\rightarrow 0$. In
particular this shows that the cohomology classes of the cocycles
$c_l, 2\leq l\leq n+1$ attached to the gerbed tower
$F_n\rightarrow F_{n-1}...\rightarrow F_1\rightarrow E$ are
independent of the choices made to construct them.}

\medskip

{\bf Proof.}

To construct the cocycle $c_{l+2}$ we pick elements
$c_{i_1...i_{l+2}}$ which represents $2$-arrows of an object
${u_{i_1..i_{l+1}}}^{i_{l+2}}$ that we lift to $2$-arrows
${c_{i_1..i_{l+2}}}^*$ of the object $u_{i_1..i_{l+2}}$  of
${F_{ll-1}}_{{{u_{i_1..i_{l+1}}}}^{i_{l+2}}}$. The representant
$c_{i_1..i_{l+3}}$ of $[c_{l+2}]$ are the Cech boundary of the
family ${c_{i_1..i_{l+2}}}^*$ acting on $u_{i_1..i_{l+2}}$. This
is by definition the construction of the connecting morphism
$H^{l+1}(E,L_l)\rightarrow H^{l+2}(E,L_{l+1})$ of the exact
sequence $0\rightarrow L_{l+1}\rightarrow L'_{l+1}\rightarrow
L_l\rightarrow 0$ $\bullet$

\medskip

{\bf Definition 4.2.4.}

A gerbed tower $f_n=F_n\rightarrow F_{n-1}...F_1\rightarrow E$ is
trivial if and only if there exists a gerbed tower
$f_{n-1}=F'_{n-1}\rightarrow F'_{n-2}\rightarrow ..F'_1\rightarrow
E$ and an element   $p\in\{1,..,n-2\}$ such that $F'_l$ is $F_l$
if $l\leq p$, $F'_p$ is a sub $2$-category of $F_{p+1}$, for each
$2$-arrow $h$ of $F_{p-1}$, the image of the arrow $h^*$ of
${F'_p}$ by $p_{p+1}$ is an arrow $u^*$, of $F_p$, defined by the
axiom $5$,  defined by
  a $2$-arrow of $u$ of $F_p$. For $l>p$, the
  $2$-category $F'_l$ is a subcategory of $F_{l+1}$.  The gerbed tower $f_{n-1}$
is called a trivialization of $f_n$ $\bullet$

\medskip

{\bf Proposition 4.2.5.}

{\it The class $[c_{n+1}]$ of a trivial gerbed tower
$F_n\rightarrow..\rightarrow F_1\rightarrow E$ is zero.}

\medskip

{\bf Proof.}

Let $f_n=F_n\rightarrow F_{n-1}\rightarrow... F_1\rightarrow E$ be
a trivial gerbed tower, and $F'_{n-1}\rightarrow
F'_{n-2}...\rightarrow F'_1\rightarrow E$ a trivialization of
$f_n$. Suppose that the integer $p$ of the definition above is
$n-2$, this means that if $l\leq n-2$, and $F'_l$ is $F_l$. We
denote by $(L_1,...,L_n)$ the band of the gerbed tower
$F_n\rightarrow...F_1\rightarrow E$, and by $L"_n$ the  sheaf such
that the group $Aut({F'_{n-1n-2}}_{h_{n-1}})$ of automorphisms of
a $2$-arrow $h_{n-1}$ of the object $U_{n-1}$ of $F'_{n-1}$ which
project by $p'_1..p'_{n-1}$ to elements of
$L_{n-2}(p'_1..p'_{n-1}(U_{n-1}))$ is
$L"_n(p_1..p'_{n-1}(U_{n-1}))$. We have the commutative:

$$
\matrix{ 0\longrightarrow  L_n & \longrightarrow  L'_n
&\longrightarrow L_{n-1} &\longrightarrow 0\cr \downarrow
&\downarrow &\downarrow\cr 0\longrightarrow  L_n &\longrightarrow
 L"_n &\longrightarrow  L_{n-2}& \longrightarrow 0}
$$

The map of $L'_n\rightarrow L"_n$ is defined by the restriction of
the morphisms $u^*$, where $u\in L_{n-1}$, and the map
$L_{n-1}\rightarrow L_{n-2}$ is zero. This exact sequence gives
rise to the commutative diagram:

$$
\matrix{H^n(E,L_{n-1}) &{{\longrightarrow}} & H^n(E,L_{n-2})\cr
\downarrow & \downarrow \cr H^{n+1}(E,L_{n}) & {{\longrightarrow}}
& H^{n+1}(E,L_{n})}
$$

Since the map $H^n(E,L_{n-1}){{\longrightarrow}}H^n(E,L_{n-2})$ is
zero, and the map $H^{n+1}(E,L_{n}) {{\longrightarrow}}
H^{n+1}(E,L_{n})$ is the identity, we deduce that the class
$[c_{n+1}]$ of the classifying cocycle of the gerbed tower
$F_n\rightarrow F_{n-1}...\rightarrow F_1\rightarrow E$ is zero.

If $p$ is not $n-1$, the last argument show that $[c_{p+1}]$ is
zero. This implies that $[c_l]=0$ for $l\geq p+1$ $\bullet$

\bigskip

{\bf Theorem 4.2.6.}

{\it Let $F_n\rightarrow F_{n-1}...F_1\rightarrow E$ be a gerbed
tower bounded by the family of sheaves $(L_1,...,L_{n})$  whose
classifying cocycles are $(c_2,...,c_{n+1})$. Consider  an exact
sequence of sheaves $0\rightarrow L_{n+1}\rightarrow
L'_{n+1}\rightarrow L_n\rightarrow 0$. Then there exists a gerbed
tower $F_{n+1}\rightarrow F_n...\rightarrow F_1\rightarrow E$
whose classifying cocycles are $c_2,...,c_{n+1},c_{n+2}$, where
$c_{n+2}$ is a $n+1$-cocycle whose cohomology class is the image
of the class of $c_{n+1}$ by the connecting map
$H^{n+1}(E,L_n)\rightarrow H^{n+2}(E,L_{n+1})$ defined by the
previous exact sequence.}

\medskip

{\bf Proof.}

Let $U_n$ be an object of $F_n$, and $U_{n-1}$ its image by the
projection map $p_n:F_n\rightarrow F_{n-1}$. The topos $U_n$ is a
$L_n$-torsor defined over an object of the topos $U_{n-1}$ since
${F_{nn-1}}_{U_{n-1}}$ is a $L_n$-gerbe defined on $U_{n-1}$. An
object $V_n$ of ${F_{nn-1}}_{U_{n-1}}$ is a $L_n$-torsor defined
over an object $U'_{n-1}$ of the topos $U_{n-1}$. This torsor is
defined by a trivialization $(V_i)_{i\in I}$, and coordinate
changes $u_{ij}:V_i\times_{U_{n-1}}V_j\rightarrow L_n$. This
coordinate changes define a principal $L_n$-torsor over $U'_{n-1}$
by gluing $V_i\times L_n$. Without restricting the generality, we
can suppose that the objects of ${F_{nn-1}}_{U_{n-1}}$ are
$L_n$-principal torsors.  We define the fiber ${F_{n+1}}_{U_n}$ to
be the gerbe bounded by $L_{n+1}$ defined on $U_n$ which
represents the obstruction to lifts the $L_n$-torsor
$p_{U_n}:U_n\rightarrow U'_{n-1}$ to a $L'_{n+1}$-torsor whose
quotient by $L_{n+1}$ is the previous $L_n$-torsor
$p_{U_n}:U_n\rightarrow U'_{n-1}$.

The objects of the category of morphisms
$Hom_{F_{n+1}}(U_{n+1},U'_{n+1})$ between the objects $U_{n+1}$ of
${F_{n+1n}}_{U_n}$ and $U'_{n+1}$ of ${F_{n+1n}}_{U'_n}$ are the
$2$-arrows between $U_n$ and $U'_n$ of $F_n$, the $2$-arrows are
the morphisms of torsors ${u^2}_{n+1}$ such that there exists a
$2$-arrow $u^1_n:U_n\rightarrow U'_n$ such that the following
diagram is commutative:

$$
\matrix{{{U}_{n+1}}& {\buildrel{{u^2_n}}\over{\longrightarrow}}&
{{U'}_{n+1}}\cr \downarrow p_{n+1} && \downarrow p_{n+1}\cr {U}_n
& {\buildrel{u^1_n}\over{\longrightarrow}} & {U'}_n}
$$

We show now that $F_{n+1}\rightarrow F_n...\rightarrow
F_1\rightarrow E$ is a gerbed tower:

Let $U_{n}$ be an object of $F_{n}$, and $U$ the object
$p_1..p_{n}(U_{n})$ of $E$. For each map $f:V\rightarrow U$, we
can define the restriction $r^n_{V,U}(f):{F_n}_U\rightarrow
{F_n}_V$. The restriction $r^{n+1}_{V,U}$ is defined on
${F_{n+1}}_{U_n}$ by the pull-back of $U_{n+1}$ by the arrow
$r^n_{V,U}(f)$.

The definition of $F_{n+1}$ implies that for every object $U_n$ of
$F_n$, the gerbe ${F_{n+1}}_{U_n}$ is bounded by $L_{n+1}$.

For every $2$-arrow $u^2_n:h^1_n\rightarrow h^2_n$, (recall that
$u^2_n$ is a morphism between the topoi ${U_1}_n$ and ${U_2}_n$)
of $F_n$, the functor ${u^2_n}^*$ is defined as follows: Without
restricting the generality, we can suppose that ${U_1}_n$ is the
trivial torsor to $V_1\times L_n$ and ${U_2}_n$ the trivial torsor
$V_2\times L_n$; $u^2_n$ is then a morphism of $L_n$-torsors,
${u^2_n}^*$ is a morphism such that the following diagram is
commutative:

$$
\matrix{V_1\times L'_{n+1}&
{\buildrel{{u^2_n}^*}\over{\longrightarrow}}& V_2\times
L'_{n+1}\cr \downarrow  && \downarrow \cr V_1\times L_n &
{\buildrel{u^2_n}\over{\longrightarrow}} & V_2\times L_n}
$$

This shows that $F_{n+1}\rightarrow...F_1\rightarrow E$ is a
gerbed tower.

The classifying cocycle of this gerbed tower is constructed by
considering the automorphism $c_{i_1..i_{n+2}}$ of the object
${u_{i_1..i_{n+2}}}$ of $F_{n+1}$, that we suppose to be
isomorphic to a trivial torsor $V_{i_1..i_{n+2}}\times L_{n}$, the
morphism $c_{i_1..i_{n+2}}$ can be lifted to an element
${c_{i_1..i_{n+2}}}^*$ of $V_{i_1..i_{n+2}}\times L'_{n+1}$. The
Cech boundary of ${c_{i_1..i_{n+2}}}^*$ is the classifying cocycle
of the gerbed tower. The cohomology class of this cocycle is the
image of the cohomology class of $c_{n+1}$ by the connecting
morphism $H^{n+1}(E,L_n)\rightarrow H^{n+2}(E,L_{n+1})$ $\bullet$

\medskip
\bigskip

\section{ Spectral sequences and gerbed towers.}

\bigskip

The goal of this part is to apply spectral sequences to study
commutative gerbed towers.

Let $E(L_1,..,L_n,...)$ be an $\infty$-gerbed tower, where
$(L_n)_{n\in {\N}}$ is a family of commutative sheaves defined on
$E$. We suppose that the topology of $E$ is defined by the
covering family $(X_i\rightarrow X)_{i\in I}$. We define
$L=\oplus_{i\geq 1} L_i$, and denote by $(C^*(X,L),d)$ the complex
of Cech $L$-chains defined on $E$. We can endow this chain complex
with the following filtration:
$$
V_p=C(X,{\oplus L_i}_{i\geq p}),
$$
and with the graduation
$$
V^p=C^p(X,L)
$$
We will calculate the terms associated to the spectral sequence
associated to this graduation.

Denote by $ Z^p_r=\{x\in V_p: d(x)\in V_{p+r}\},$
$B^p_r=d(V_{p-r})\cap V_p,$
   and $E^p_r= {Z^p_r\over
{Z^{p+1}_{r-1}+B^p_{r-1}}}$.

We suppose in the sequel that $r\geq 1$.

\medskip

{\bf Determination of $Z_r^p$.}

\medskip

Let $x$ be an element of $V_p$, $d(x)$ is an element of $Z_r^p$ if
and only if $d(x)$ is an element of $V_{p+r}$. We can write
$x=x_{i_p}+..+x_{i_n}$, where $x_{i_l}$ is the homogeneous
component of $x$ which takes value in $L_l$, $d(x_l)$ is an
element of $V_{p+r}$ if and only if $d(x_l)=0$ if $l\leq p+r$. We
deduce that $x$ is an element of $Z_r^p$ if and only if its
components $x_{i_j}$, such that $j<p+r$ are cocycles.

\medskip

{\bf Determination of $B^p_r$.}

\medskip

Let $x$ be an element of $V_{p-r}$, and $x_l$ its component which
takes values in $L_l$. The image by $d$ of $x_l$ is an element of
$B_r^p$ if and only if $d(x_l)$ is an element of $V_p$. This
equivalent to saying that $d(x_l)$ is zero, or $l\geq p$. This
implies that $B_r^p=d(V_p)$.

\medskip

{\bf Determination of $E_r^p$.}

\medskip

We have $Z^p_r=Z^{p+1}_{r-1}\oplus Z(V_p)\cap C(X, L_p)$. We
deduce that $E_r^p=H(X,L_p)$.

\bigskip

Now we set $Z^{pq}_r=Z^p_r\cap V^{p+q}$, $B_r^{pq}=B_r^p\cap
V^{p+q}$ and $E^{pq}_r={Z^{pq}_r\over
B^{pq}_{r-1}+Z^{p+1,q-1}_{r-1}}$.

\medskip

{\bf Determination of $Z_r^{pq}$.}

\medskip

Let $x$ be an element of $Z_r^p$, one of its homogeneous
components $x_l$ which takes values in $L_l$, is an element of
$Z_r^{pq}$, if $x_l$ is a $V_p$ $p+q$-chain, and $d(x_l)$ is an
element of $V_{p+r}$. We have seen that $d(x_l)$ is an element of
$Z^p_r$ if and only if $x_l$ is a cocycle or $r\geq p+r$. We
deduce from this fact that $Z^{pq}_r=C^{p+q}(X,V_{p+r})\oplus
Z^{p+q}(X,L_1\oplus..\oplus L_{p+r-1})$.

\medskip

{\bf Determination of $B_r^{pq}$.}

\medskip

Let $x$ be and element of $B_r^{pq}$, one of its homogeneous
components, $x_l$ which takes value in $L_l$ is an element of
$B_r^{pq}$ if and only if it is a $p+q$-chain, and there exists an
element $y$ in $V_{p-r}$ such that $d(y)=x_l$. We deduce that
$B^{pq}_r=d(C^{p+q-1}(X,V_p))$.

\medskip

{\bf Determination of $E_r^{pq}$.}

\medskip

The vector space $Z_r^{pq}$ is the summand of $Z^{p+1,q-1}_{r-1}$
and  $Z^{p+q}(X,L_p)$. We deduce that $E^{pq}_r=H^{p+q}(X,L_p)$.

\medskip

Now, we will denote by $Z^p_{\infty}$, the set of cocycles
contained in $V_p$, by $B^p_{\infty}$ the set of boundaries
contained in $V_p$, and by
$E^p_{\infty}={Z^p_{\infty}\over{Z^{p+1}_{\infty}+B^p_{\infty}}}$.
We remark that $E^p_{\infty}=H(X,L_p)$.

The following proposition can be deduced from [19] p. 84 Theorem
4.6.1.

\medskip

{\bf Proposition 4.0.1.}

{\it Suppose that there exists an integer $n\geq r$ such that
$H^{p+q}(X,L_p)=E^{pq}_r=0$ for $p\neq 0,n$ and an integer $s$
such that $L_n=0$ if $n>s$, then we have the following exact
sequence $$ ...\rightarrow H^{i}(X,L_n)\rightarrow
H^i(X,L)\rightarrow H^i(X,L_1)\rightarrow
H^{i+1}(X,L_n)\rightarrow H^{i+1}(X,L)\rightarrow...
$$
}
\bigskip

\section{ Application of gerbed towers to affine manifolds.}

\bigskip

 An affine manifold $(N,\nabla_N)$, is a differentiable manifold
 $N$, endowed with a connection $\nabla_N$ whose curvature and
 torsion forms vanish identically. The connection $\nabla_N$ defines
 on $N$ an atlas whose coordinate changes are affine transformations.
 Auslander has conjectured that the fundamental group of a compact
 and complete affine manifold is polycyclic. Let $(N,\nabla_N)$ and
 $(N',\nabla_{N'})$ be two affine manifolds of respective dimension $n$ and $n'$ whose affine structures are defined
 by the respective atlases $(U_i,u_i)$, and $(U'_j,u'_j)$. An affine map $f:(N,\nabla_N)\rightarrow
 (N',\nabla_{N'})$ is a differentiable map $f:N\rightarrow N'$
 such that $u'_i\circ f_{\mid U_i}\circ {u_i}^{-1}$
 is a restriction of an affine map from ${\R}^n$ to ${\R}^{n'}$.
 Suppose that $(N,\nabla_N)$ and $(N',\nabla_{N'})$ are complete
 and compact. It is shown in Tsemo [26] that in this case,  there exists a
compact and complete affine manifold $(N_1,\nabla_{N_1})$ of
dimension $n'$, and an affine submersion
$f_1:(N,\nabla_N)\rightarrow (N_1,\nabla_{N_1})$. Ehresman has
shown that submersions between compact  manifolds are locally
trivial differentiable fibrations. The typical fiber $F$ of the
fibration $N\rightarrow N_1$ inherits from $N$ complete affine
structures.  The homotopy exact sequence of this fibration gives
rise to the sequence:

$$
1\longrightarrow
\pi_1(F)\rightarrow\pi_1(N)\rightarrow\pi_1(N_1)\rightarrow 1
$$
If the fundamental groups of $\pi_1(F)$ and $\pi_1(N_1)$ are
polycyclic, then $\pi_1(N)$ is also polycyclic. This has motivated
the following conjecture:

\medskip

{\bf Conjecture 5.0.1.}

{\it Let $(N,\nabla_N)$ be a $n$-dimensional compact and complete
affine manifold, then there exists a complete affine structure
$(N',\nabla_{N'})$  defined on a finite galois cover $N'$ of $N$,
and a non trivial affine map $f:(N',\nabla_{N'})\rightarrow
(N_1,\nabla_{N_1})$. Non trivial means that the dimension of the
fibers of $f$ are different from zero, and $n$.}

\medskip

This conjecture implies the Auslander conjecture, and leads to the
problem of classifying sequences of affine submersions
$(N_n,\nabla_{N_n})\rightarrow..\rightarrow (N_1,\nabla_{N_1})$.
This if the last conjecture is true, will allow to know the
topology of all the compact and complete affine manifolds. The
theory of gerbed towers has been first constructed to study this
classification problem.

\medskip

\medskip

{\bf Definition 5.0.2.}

An affinely locally trivial affine fibration, whose typical fiber
is the affine manifold $(F,\nabla_F)$, is an affine map
$f:(N_1,\nabla_{N_1})\rightarrow (N,\nabla_N)$ which is the total
space of a bundle  whose fibers inherit from $(N_1,\nabla_{N_1})$,
affine structures whose holonomies is the holonomy of the affine
structure $(F,\nabla_F)$ $\bullet$

\medskip

We will restrict to  the study of sequences
$(N_n,\nabla_{N_n})\rightarrow ...(N_1,\nabla_{N_1})\rightarrow
(N,\nabla_N)$ where each map $f_p:(N_p,\nabla_{N_p})\rightarrow
(N_{p-1},\nabla_{N_{p-1}})$ is an affinely locally trivial affine
fibration.

\medskip

Let $f:(N_1,\nabla_{N_1})\rightarrow (N,\nabla_N)$ be an affinely
locally trivial affine fibration whose typical fiber is the affine
manifold $(F,\nabla_F)$. We suppose that the affine structure of
$(N_1,\nabla_{N_1})$ is complete. This implies that the affine
structure of $(N,\nabla_N)$ is complete see Tsemo [27]. We can
identify $\pi_1(N_1)$ with its image by the holonomy morphism of
$(N_1,\nabla_{N_1})$. Suppose that the dimension of $N_1$ and $N$
are respectively $n_1$ and $n$. Let $h:\pi_1(N_1)\rightarrow
Aff({\R}^{n_1})$ be the holonomy representation of
$(N_1,\nabla_{N_1})$. We can write see Tsemo [27]
${\R}^{n_1}={\R}^n\times{\R}^p$, and for each element $\gamma$ of
$\pi_1(N_1)$, $\gamma(x,y)=({L_1}_{\gamma}(x)+{l_1}_{\gamma},
{L_2}_{\gamma}(y)+{L_3}_{\gamma}(x)+{l_2}_{\gamma})$. Where
${L_1}_{\gamma}$ and ${L_2}_{\gamma}$ are respective automorphisms
of ${\R}^n$ and ${\R}^p$, $L_3:{\R}^n\rightarrow {\R}^p$ is a
linear map, and ${l_1}_{\gamma}$ and ${l_2}_{\gamma}$ are
respective elements of ${\R}^n$ and ${\R}^p$. An element $\gamma$
of $\pi_1(F)$ is an element $\gamma$ of $\pi_1(N_1)$ such that
${L_1}_{\gamma}$ is the identity of ${\R}^n$, and ${l_1}_{\gamma}$
is zero, and ${L_3}_{\gamma}$ is zero. We can identify $\pi_1(N)$
to the set of affine transformations
$({L_1}_{\gamma},{l_1}_{\gamma})$.

 Let $T_F$ be the translation group of $(F,\nabla_F)$, that is
the group of affine automorphisms of $(F,\nabla_F)$ whose elements
lift to translations of ${\R}^p$. Since the group $\pi_1(F)$ is a
normal subgroup of $\pi_1(N_1)$, the holonomy of
$(N_1,\nabla_{N_1})$ induces a representation $\pi_1(N)\rightarrow
Aff(F,\nabla_F)/T_F$, which defines a flat
$Aff(F,\nabla_F)/T_F$-bundle $p_F$ over $(N,\nabla_N)$. The
composition of the holonomy of $(N_1,\nabla_{N_1})$, and the
conjugation of $Aff(F,\nabla_F)$ defined  a flat bundle
$T_F$-bundle $p'_F$ over $(N,\nabla_N)$ see Tsemo [27]. An
isomorphism $h:e\rightarrow e'$ between a pair of  locally trivial
$(F,\nabla_F)$-affine bundles $e$ and $e'$ defined over
$(N,\nabla_N)$, is an affine map $h:e\rightarrow e'$ which is an
isomorphism of bundles which gives rise to the identity of $p_F$.

Given affine manifolds $(F,\nabla_F)$ and $(N,\nabla_N)$ and a
flat $Aff(F,\nabla_F)/T_F$-bundle $p_F$, we can define the first
extension problem as follow: study the existence and classify
affinely locally trivial affine fibrations
$f:(N_1,\nabla_{N_1})\rightarrow (N,\nabla_F)$ which give rise to
the flat bundle $p_F$.

\medskip

{\bf Proposition 5.0.3.}

{\it Let $(F,\nabla_F)$ and $(N,\nabla_N)$ be compact and complete
affine manifolds, and $p_F$ a flat $Aff(F,\nabla_F)/T_F$-bundle
defined on $N$, The bundle $p_F$ induces a flat $T_F$-bundle
$p'_F$. For each open subset $U$ of $N$, define the category
$C_F(U)$ to be the category whose objects are affinely locally
trivial $(F,\nabla_F)$-affine bundles
 which induce the restriction of $p_F$ to $U$. A map
$h:e_U\rightarrow e'_U$ between two objects of $C_F(U)$ is an
isomorphism of affine bundles which gives rise to the identity of
the restriction of $p_F$ to $U$. The correspondence defined
 on the category of open subsets of $N$, by $U\rightarrow C_F(U)$, is a gerbe whose
classifying cocycle represents the obstruction to the existence of
an affinely locally trivial affine bundle, which gives rise to
$p_F$. The gerbe $C_F(U)$ is a gerbe bounded by the sheaf of
affine sections of $p'_F$, that we denote $L_F$. }

\medskip

{\bf Proof.}

Gluing property for objects.

Let $(U_i)_{i\in I}$ be an open covering of an open subset  $U$ of
$N$, $e_i$ an object of $C_F(U_i)$, and $u_{ij}:e^i_j\rightarrow
e^j_i$ an isomorphism
 such that
 ${u_{i_1i_2}}^{i_3}{u_{i_2i_3}}^{i_1}={u_{i_1i_3}}^{i_2}$. The definition of a bundle
 implies the existence of a   bundle $e$ defined over $U$ whose restriction to $U_i$ is
 $e_i$. Since the coordinate changes $u_{ij}$ are affine
 isomorphisms between affinely locally trivial
 $(F,\nabla_F)$-affine bundles, this implies that $e$ is an
 affinely locally trivial $(F,\nabla_F)$-affine bundle.

 Gluing conditions for arrows.

 Let $e$ and $e'$ be two objects of $C_F(U)$,  the
 correspondence defined on the open subsets of $U$ by
 $V\rightarrow Hom(e_{\mid V},e'_{\mid V})$, where $e_{\mid V}$ and
 $e'_{\mid V}$ are the respective restrictions of $e$ and $e'$ to
 $V$ is a sheaf of sets, since it is a sheaf of morphisms between
 two bundles.

 \medskip

 This shows that $C_F$ is a sheaf of categories. It remains to show
 that $C$ is a gerbe.

Let $(U_i)_{i\in I}$ be an open covering of $N$ by contractible
open subsets which are domain of affine charts. Then $U_i\times
(F,\nabla_F)$ is an object of $C_F(U_i)$.

Let $U$ be an open subset of $N$, and $e$ and $e'$ a pair of
objects of $C_F(U)$. The respective restrictions $e_{\mid U_i\cap
U}$, and $e'_{\mid U_i\cap U}$ of $e$ and $e'$ to $U_i\cap U$ are
isomorphic to $U_i\cap U\times (F,\nabla_F)$.

An isomorphism $h$, of an object $e$ of $C_F(U)$, is an
isomorphism of affine bundle which gives rise to the identity of
the restriction of $p_F$ to $U$. The restriction of $h$ to
$e_{\mid U\cap U_i}$ is an isomorphism $h_i$ of the trivial bundle
$U\cap U_j\times (F,\nabla_F)$.  The fact that $h$ gives rise to
the identity of $p_F$, is equivalent to the fact that its
restriction to a fiber yields to the identity of
$Aff(F,\nabla_F)/T_F$. This implies that $h_i$ is a $T_F$ valued
affine map, and $h$ is a section of $p'_F$ $\bullet$

\medskip

{\bf Proposition 5.0.4.}

{\it Suppose that the gerbe $C_F$ is trivial, then the objects of
$C_F(N)$ are diffeomorphic manifolds.}

\medskip

{\bf Proof.}

Suppose that the gerbe $C_F$ is trivial, then the holonomy of a
global object $(N_1,\nabla_{N_1})$ is defined by a representation
$h_{\gamma}(x,y)=({L_1}_{\gamma}(x)+{l_1}_{\gamma},
{L_2}_{\gamma}(y)+{L_3}_{\gamma}(x)+{l_2}_{\gamma})$ which defines
an $(F,\nabla_F)$-bundle. The objects of $C_F(N)$ are classified
by $H^1(N,p'_F)$, the $1$-cohomology group of the sheaf of affine
sections of $p'_F$. An element of $H^1(N,p'_F)$ is defined by an
affine $C_3:{\R}^n\rightarrow {\R}^p$ which is a $1$-cocycle for
the action of $\pi_1(N)$ on $Aff({\R}^n,{\R}^p)$ defined by
$\gamma(C_3)={L_2}_{\gamma}\circ C_3\circ
{({L_1}_{\gamma},{l_1}_{\gamma})}^{-1}$. The bundle defined by the
representation
$h^t_{\gamma}(x,y)=({L_1}_{\gamma}(x)+{l_1}_{\gamma},
{L_2}_{\gamma}(y)+{L_3}_{\gamma}(x)+tC_3(({L_1}_{\gamma},{l_1}_{\gamma})(x))+{l_2}_{\gamma})$,
$0\leq t\leq 1$ defines an homotopy between the bundle $h'_1$
defined by $h^1$ and the one $h'_0$ defined by $h^0$. We deduce
that $h'_0$ and $h'_1$ are isomorphic differentiable bundles. This
implies that $N_1$ is diffeomorphic to the $F$-bundle defined by
$h^0$ $\bullet$

\bigskip

Let $(N,\nabla_F)=(F_0,\nabla_{F_0})$,
$(F_1,\nabla_{F_1})$,...,$(F_n,\nabla_{F_n})$ be affine manifolds.
We are going to define a gerbed tower which will allow us to study
the classification of sequences $(N_n,\nabla_{N_n})\rightarrow
(N_{n-1},\nabla_{N_{n-1}})\rightarrow...(N,\nabla_N)$, where
($h_l:(N_l,\nabla_{N_l})\rightarrow (N_{l-1},\nabla_{N_{l-1}})$ is
an affinely locally trivial affine bundle).

 Denote by $T_{F_l}$ the group of translations of $(F_l,\nabla_{F_l})$. We suppose defined a flat
$Aff(F_l,\nabla_{F_l})/T_{F_l}$-bundle $p_{F_l}$ over
$(F_{l-1},\nabla_{F_{l-1}})$. The bundle $p_{F_l}$ induces a gerbe
$C_l$ defined over $F_{l-1}$ (see proposition 4.3).

\medskip

{\bf Definition 5.0.5.}

We define  $L_n\rightarrow L_{n-1}...\rightarrow L_1\rightarrow
L_0$ to be the sequence of $2$-categories such that
$L_1\rightarrow L_0$ is $C_{1}$, supposed defined $L_p$, an object
of $L_p$ is an affinely locally trivial
$(F_p,\nabla_{F_p})$-bundle. A $1$-arrow between a pair of objects
$e_p$ and $e'_p$ of $L_p$ is an affine map between their base
space. A $2$-arrow is an isomorphism of affinely locally trivial
bundle which cover a $1$-arrow. We denote by $p_l:L_l\rightarrow
L_{l-1}$ the projection.

An object of $L_{p+1}$ is an affinely locally trivial
$(F_{p+1},\nabla_{F_{p+1}})$-bundle $e_{p+1}$ defined over an open
subset $U_{p}$ of  an object $e_p$ of $L_p$ such that the
restriction to each fiber of $U_p\rightarrow p_p(e_{p})$ of the
$Aff(F_{p+1},\nabla_{F_{p+1}})/T_{F_{p+1}}$-bundle induced, is the
restriction of the
$Aff(F_{p+1},\nabla_{F_{p+1}})/T_{F_{p+1}}$-bundle induced is
$p_{F_{p+1}}$. A $1$-arrow $h^1_{p+1}:e_{p+1}\rightarrow e'_{p+1}$
between a pair of objects $e_{p+1}$ and $e'_{p+1}$ of $L_{p+1}$ is
a affine map between their respective base spaces $e_p$ and $e'_p$
induced by a $2$-arrow of $L_p$. A $2$-arrow between $e_{p+1}$ and
$e'_{p+1}$ is an isomorphism of affinely locally trivial
$(F_{p+1},\nabla_{F_{p+1}})$-bundles which covers a $1$-arrow. We
suppose that the $2$-arrows depend only of the open subset
$p_1..p_{p+1}(e_{p+1})$ and $p_1..p_{p+1}(e'_{p+1})$. Let
$h^2_{p}$ be a $2$-arrow of $L_p$, we define ${h^2_{p}}^*$ to be a
$2$-arrow of $L_{p+1}$ which cover $h^2_{p}$ $\bullet$

\medskip

{\bf Proposition 5.0.6.}

{\it The sequence $L_n\rightarrow L_{n-1}...\rightarrow
L_1\rightarrow L_0$ that we have just defined is a gerbed tower.}

\medskip

{\bf Proof.}

The fibered category $L_1\rightarrow L_0$ is a gerbe as shows
proposition 4.4. Let $V$ be an open subset of $N$, and $U$ an open
subset of $V$, the restriction functor $r_{U,V}:{L_p}_U\rightarrow
{L_p}_V$ is defined by the restriction of bundles.

Let $U$ be an open subset of $N$, and $e_l$ be an object of
${L_l}_{U}$, the band of ${L_{l+1}}_{e_l}$ is the sheaf of
sections of the $T_{F_{l+1}}$-bundle $p'_{e_l}$ defined on $e_l$
induced by $p_l$, this   sheaf does not depend of the objects
chosen in the fibre ${L_{l+1}}_{e_l}$, since we have supposed that
the $2$-arrows depend only of $N$.

Consider $2$-morphisms $u_p:e_p\rightarrow e'_p$, and
$u'_p:e'_p\rightarrow e"_p$ of $L_p$. We have defined in the
paragraph above the proposition a morphism ${u_p}^*$. These
morphisms satisfy ${u'_p}^*{u_p}^*=c(u_p,u'_p)({u'_pu_p})^*$,
where $c(u_p,u'_p)$ is an automorphism of an object of the gerbe
${L_{p+1}}_{e"_p}$ induced by the band $\bullet$

\bigskip

\section{ Interpretation of the integral cohomology of a
manifold.}

\bigskip

Characteristic classes have been used by many mathematicians to
study geometric objects. On this purpose, we have to give a
geometric interpretation of the group $H^n(N,{\Z})$. This is what
we propose to do in this part.

It is a well-known fact that the group  $H^2(N,{\Z})$ is the set
of equivalence classes of complex line bundles over $N$. Brylinski
has defined an equivalence between the space of equivalence
classes of complex line gerbes and $H^3(N,{\Z})$.

Consider ${\C}^*_N$ the sheaf of differentiable
${\C}-\{0\}={\C}^*$-functions defined on $N$. We say that a class
$[c_n]$ of $H^n(N,{\C}^*_N)$ is geometric, if and only if there
exists an $(n-1)$-gerbed tower which classifying cocycle is $c_n$.
A sufficient condition for a class $c_n$ to be geometric is the
following: there exists a classifying cocycle $c_{n-1}$ of a
commutative $(n-2)$-gerbed tower $E(L_1,..,L_{n-2})$ which  is an
element of $H^{n-1}(N,L_{n-2})$, an exact sequence of sheaves
$0\rightarrow {\C}^*_N\rightarrow L\rightarrow L_{n-2}\rightarrow
0$ such that $[c_{n}]$ is the image of $[c_{n-1}]$ by the boundary
map $\delta:H^{n-1}(N,L_{n-2})\rightarrow H^n(N,{\C}^*)$.

We have the exact sequence
$$
0\rightarrow {\Z}{\buildrel  i\over{\rightarrow}} {\C}_N{\buildrel
{exp } \over{\rightarrow}}{\C}^*{_N}\rightarrow 0.
$$
where $i$ is the canonical injection, and $exp$ the exponential
map. It results from this sequence an isomorphism between
$H^n(N,{\C}^*_N)$ and $H^{n+1}(N,{\Z})$. An element of
$H^{n+1}(N,{\Z})$ will be said geometric if and only if it is the
image of an element of $H^n(N,{\C}^*_N)$ which is the classifying
cocycle of a $(n-1)$-gerbed tower by the preceding isomorphism.

\bigskip

We have the following result:

\bigskip

{\bf Theorem 6.0.1.}

{\it Let $N$ be a differentiable manifold, then each geometric
class of $H^{n+2}(N,{\Z})$ is the classifying cocycle of a
$n$-gerbed tower defined on $N$.}

\bigskip

\section{  $n$-categories, and sheaves of $n-$categories.}

\bigskip

In this part, we will define a notion of sheaf of $n$-categories
over a topos  $N$.

\medskip

{\bf Definition 7.0.1.}

A $0$-pseudo-category is a set, a $1$-pseudo-category  $C_1$, is a
category.

Suppose defined the notion of $n$-pseudo-category.

An $(n+1)$-pseudo-category  $C_{n+1}$, is defined by, a class of
objects $Ob(C_{n+1})$, for each objects $x$, and $y$, the
$n$-pseudo-category  of morphisms $Hom(x,y)$. For each objects
$u_1,u_2$ and $u_3$ of $C_{n+1}$, there exists a composition
$n$-functor:

$$
Hom(u_2,u_3)\times Hom(u_1,u_2)\longrightarrow Hom(u_1,u_3)
$$
 We
suppose the existence of an object $1_x$ of $Hom(x,x)$, such that
for each arrow $h:x\rightarrow y$, $h\circ 1_x$ is isomorphic to
$h$, and for each arrow $h':y'\rightarrow x$, $1_xh'$ is
isomorphic to $h'$.

An isomorphism between the objects $x$ and $y$ of an
$n+1-$pseudo-category is a map $f:x\rightarrow y$, such that there
exists $h:y\rightarrow x$ such that $hf$ is isomorphic to $1_x$,
and $fh$ to $1_y$.

A functor between two $n$-pseudo-categories  $C_n$ and $C'_n$, is
defined as follows:

(i) A map $F:Ob(C_n)\rightarrow Ob(C'_n)$, and for each arrow
$f:x\rightarrow y$, a morphism $F(f):F(x)\rightarrow F(y)$ such
that $F(f\circ f')$ is isomorphic to $F(f)\circ F(f')$.

(ii) A natural transformation between two functors $F$ and $F'$,
is defined by a family of maps $u_x:F(x)\rightarrow F'(x)$ such
that for each map $f:x\rightarrow y$, $u_yF(f)$ is isomorphic to
$F'(f)u_x$ At this stage, we do not precise the gluing datas.

\medskip

\medskip

{\bf Definition 7.0.2.}

A $0$-sheaf of sets defined on $N$, will be a sheaf of sets.
Suppose defined the notion of sheaves of $n-1$-pseudo-categories.
A sheaf of $n$-pseudo-categories, will be defined by the following
data: for each pair of objects $U$ and $V$ of $N$, and a map
$h:U\rightarrow V$,  a restriction functor
${r^{C_n}}_{U,V}(h):C_n(V)\rightarrow C_n(U)$ such that  for each
triple of objects $U_1,U_2$ and $U_3$, there exists an isomorphism
$c(U_1,U_2,U_3)$ between
${r^{C^n}}_{U_1,U_2}(h){r^{C^n}}_{U_2,U_3}(f)$ and
${r^{C^n}}_{U_1,U_3}(fh)$.

Gluing condition for objects:

 Let $(U_i)_{i\in I}$ be an open cover of the object $U$ of $N$,
 $x_i$ an object of $C_n(U_i)$. If there exists, a family of
maps $u_{jl}:{r^{C^n}}_{U_l\times_NU_j,U_l}(x_l)\rightarrow
{r^{C^n}}_{U_l\times_NU_j,U_l}(x_j)$, a sheaf of
$n-1$-pseudo-categories $C_{n-1}$ defined on $U$, such that the
restrictions map ${r^{C_{n-1}}}_{{U_i\times_NU_j,U_i}}$ are
$u_{ji}$, then there exists an object $x$ of $C_n(U)$ such that
the restriction of $x$ to $U_i$ is $x_i$.

Gluing conditions for arrows:

For each $x$, and $y$ in $C_n(U)$, the map defined on sub-objects
of $U$ by $V\rightarrow Hom(x_{\mid V},y_{\mid V})$ is a sheaf of
$(n-1)$-categories $\bullet$

\medskip

 We denote ${\N}_n$ the pseudo-category whose objects are the elements  of the set
$\{1,...,n\}$, $Hom(j_1,j_2)$ has one element if $j_1$ inferior to
$j_2$, if not it is empty.  We endow it with the topology such
that the covering family of $l$ are the integers inferior to $l$.

\bigskip

 {\bf Definition 7.0.3.}

An $n$-category, is a $n$-pseudo-category, such that for any
objects $x_0,...,x_{n},..$ of $C$, for each family of  maps
$u_{ij}:x_j\rightarrow x_i$ such that
$u_{i_1i_2}u_{i_2i_3}=u_{i_1i_3}$, the map defined on ${\N}$,
$i\rightarrow Hom(x_i,x_0)$, is a sheaf of $n-1$-pseudo-categories
whose restrictions functors ${u_{ij}}^*:Hom(x_j,x_0)\rightarrow
Hom(x_i,x_0)$ are defined  by: $h\rightarrow hu_{ji}$ $\bullet$

\bigskip

\end{document}